\documentclass[reqno,11pt]{amsart}
\usepackage{amsmath,amssymb,tabularx,setspace,color,multirow,graphics,graphicx}
\usepackage[margin=3.5cm]{geometry}
\newtheorem{theorem}{Theorem}
\usepackage{adjustbox}

\newtheorem{definition}{Definition}
\newtheorem{remark}{Remark}
\newtheorem{corollary}{Corollary}
\makeatletter
\renewcommand\section{\@startsection {section}{1}{\z@}%
                                   {-3.5ex \@plus -1ex \@minus -.2ex}%
                                   {2.3ex \@plus.2ex}%
                                   {\normalfont\large\bfseries}}
\makeatother

\newtheorem{Corollary}{Corollary}

\newtheorem{Theorem}{Theorem}

\newtheorem{Example}{Example}

\begin{document}
\doublespace
\title[]{Weighted cumulative residual Entropy Generating Function and its properties}
\author[]%
{S\lowercase{mitha} S.$^{\lowercase{a}}$, S\lowercase{udheesh} K. K\lowercase{attumannil}$^{\lowercase{b}}$  \lowercase{and}   S\lowercase{reedevi} E. P.$^{\lowercase{c}}$  \\
 $^{\lowercase{a}}$K E C\lowercase{ollege} M\lowercase{annanam,} K\lowercase{erala}, I\lowercase{ndia},\\
 $^{\lowercase{b}}$I\lowercase{ndian} S\lowercase{tatistical} I\lowercase{nstitute},
  C\lowercase{hennai}, I\lowercase{ndia,}\\
$^{\lowercase{c}}$C\lowercase{ochin} U\lowercase{niversity of} S\lowercase{cience and} T\lowercase{echnology}, K\lowercase{ochi}, I\lowercase{ndia.}}
\maketitle
\vspace{-0.2in}
\begin{abstract}
The study on the generating function approach to entropy become popular as it generates several well-known entropy measures discussed in the literature. In this work, we define the weighted cumulative residual entropy generating function (WCREGF) and study its properties.  We then introduce the dynamic weighted cumulative residual entropy generating function (DWCREGF).  It is shown that the DWCREGF determines the distribution uniquely. We study some characterization results using the relationship between the DWCREGF and the hazard rate and/or the mean residual life function. Using a characterization based on DWCREGF, we develop a new goodness fit test for Rayleigh distribution. A Monte Carlo simulation study is conducted to evaluate the proposed test. Finally, the test is illustrated using two real data sets.
\\\noindent Keywords: Entropy; Entropy generating function; Weighted cumulative residual entropy;  Rayleigh distribution; U-statistics.
\end{abstract}
\vspace{-0.1in}
\section{Introduction}
The concept of entropy was introduced first by Shannon (1948) for measuring the uncertainty associated with a random variable.   Let $X$ be a continuous non-negative random variable having distribution function $F(.)$,  probability density function $f(.)$ and survival function $\bar{F}(x)=P(X>x)$.    The  Shannon entropy of $X$ is defined as
\begin{equation}\label{shannon}
    H(X)=-\int_{0}^{\infty}\log f (x)f(x)dx=E(-\log f(X)),
\end{equation}
where ``log'' denotes the natural logarithm.
Different entropy measures have been introduced in the literature which is suitable for some specific situations. The widely used measure of entropy is the cumulative residual entropy (CRE)  is given by (Rao et al. 2004)
 $$ \mathcal{E}(X)=-\int_{0}^{\infty}\bar{F}(x)\log \bar{F}(x)dx.$$

Di Crescenzo and Longobardi (2009) introduced the cumulative entropy (CE) for estimating the uncertainty in the past lifetime of a product/subject and is given by
$$\mathcal{CE}(X)=-\int_{0}^{\infty}F(x)\log F(x)dx.$$
The weighted versions of  $\mathcal{E}(X)$ and $\mathcal{CE}(X)$ have been studied in the literature as well.  These are given by (Mirali et al., 2016)
$$\mathcal{E}^{w}(X)=-\int_{0}^{\infty}x\bar{F}(x)\log \bar{F}(x) dx$$
 and (Mirali and Baratpour, 2017)
$$ \mathcal{CE}^{w}(X)=-\int_{0}^{\infty}xF(x)\log F(x)dx.$$
For some recent development in this area, we refer to  Sudheesh et al. (2022), Balakrishnan et al. (2022) and Chakraborty and Pradhan (2023). Among these, Sudheesh et al. (2022) defined a generalized cumulative residual entropy and studied its properties.  They show that cumulative residual entropy, weighted cumulative residual entropy and their weighted versions are special cases of the proposed measure.

The moment generating function (m.g.f) of a probability distribution is a convenient tool for evaluating mean, variance and other moments of a probability distribution. The successive derivative of the m.g.f at a point zero gives the successive moments of the probability distribution provided these moments exist. In information theory, generating functions have been defined for probability density functions to determine information quantities such as Shannon information, extropy and Kullback-Leibler divergence. Golomb (1966) introduced the entropy generating of a probability distribution given by
$$
B(s)=\int f^{s}(x) d x,\, s>0 .
$$
The first derivative of $B(s)$ at $s=1$, gives the negative of  Shannon's entropy in (\ref{shannon}).
%

In recent times, there has been a substantial focus on information generating functions. Clark (2020) proposed an information-generating function applicable to point processes to compute statistics associated with entropy and relative entropy.   In a series of papers, Kharazmi and Balakrishnan (2021a, 2021b, 2021c) introduced measures such as cumulative residual information generating and relative cumulative residual information generating, along with cumulative residual Fisher information and relative cumulative residual Fisher information measures, and conducted a study on their properties. Specifically, Kharazmi and Balakrishnan (2021a) introduced two novel divergence measures and demonstrated that notable information divergence metrics, including Jensen-Shannon, Jensen-extropy, and Jensen-Taneja, are all specific instances of these measures.
 To explore information-generating functions linked to maximum and minimum ranked sets, as well as record values and their properties, refer to the works of Zamani et al. (2022) and Kharazmi et al. (2021).
 Saha and Kayal (2023) put forth general weighted information and relative information-generating functions and conducted an in-depth study of their properties. Capaldo et al.(2023) introduced the cumulative information generating function and its suitable distortions-based extensions. Recently, Smitha et al. (2023) proposed a cumulative residual entropy generating function and is  given by
\begin{equation}\label{cresg}
   C_s(F)= \int_{0}^{\infty}(\bar{F}(x))^{s} d x ,\,\,s>0.
\end{equation}
They also proposed the dynamic version of the $ C_s(F)$. Motivated by these recent developments we study the properties of weighted version of $C_s(F).$


The rest of the paper is organized as follows. In Section 2, we define the weighted cumulative residual entropy generating function (WCREGF) and study some properties.  In Section 3, we introduced the dynamic weighted cumulative residual entropy generating function  (DWCREGF).  It is shown that the DWCREGF determines the distribution uniquely. Also, we establish some characterization results using the relationship of DWCREGF with hazard rate and mean residual life function. In Section $4$, we propose a test for testing the Rayleigh distribution using the property of  DCWREGF.  We conduct a Monte Carlo Simulation to study the finite sample performance of the proposed test. Some concluding remarks along with some open problems are given in Section 5.

\section{Weighted cumulative residual entropy generating function}
The survival function is more useful than the probability density function in lifetime studies. This motivates us to define WCREGF and study its properties. Next, we define WCREGF.
\begin{definition}
  Let $X$ be a continuous non-negative random variable having survival function $\bar{F}(x)$, then WCREGF, denoted by $C_{s}(W,F)$  is defined as
  \begin{equation}
C_{s}(W,F)=\int_{0}^{\infty}{x}(\bar{F}(x))^{s} d x,\,\,s>0.
\end{equation}
\end{definition}
Differentiating $C_{s}(W,F)$ with respect to $s$, at $s=1$ and putting a negative sign, we obtain (prime denotes the derivative)
$$
\left.C_{s}^{\prime}(W,F)\right|_{s=1}=-\int_{0}^{\infty}{x}\bar{F}(x) \log \bar{F}(x) d x,
$$
The above measure is the weighted cumulative residual entropy introduced by Mirali et al. (2016). \\
In the following example, we show that two distributions can possess the same $C_{s}(F)$ but different $C_{s}(W,F)$, which shows the importance of the weighted generating function approach to residual entropy.
\begin{Example}
 Suppose $C_{s}(F)$ and $C_{s}(G)$ are the CREGF of the random variables $X$ and $Y$, respectively. Suppose the  density functions of $X$ and $Y$ are:
 $$
 f_{X}(x)=\frac{1}{2} ,\,0\leq x\leq 2
 $$
 and
 $$
 f_{Y}(y)=\frac{1}{2} ,\,2\leq y\leq4.
 $$
 Using (\ref{cresg}), we obtain
$$
C_{s}(X)=\frac{2}{s+1}\quad\text{and}\quad C_{s}(Y)=\frac{2}{s+1}.
$$
We also find
$$
C_{s}(W,F)=\frac{4}{(s+1)(s+2)}\quad\text{and}\quad C_{s}(W,F)=\frac{4(s+3)}{(s+1)(s+2)}.
$$
We can see that
$$
C_{s}(F)=C_{s}(G).
$$
However, the WCREGF of $X$ and $Y$ are not identical.  Even though $C_{s}(F)=C_{s}(G)$, the WCREGF of $X$ is always smaller than that of $Y$ for all values of $s>0$.
\end{Example}
In Table 1, we derive the expression for $C_{s}(W,F)$ for some well-known lifetime distributions.
\begin{table}[h]
  \centering
   \caption{Expression of $C_{s}(W,F)$ for some well-known distributions}
\begin{tabular}{cccccccc}
\hline
 Distribution & $\bar{F}(x)$ & $C_s(W,F)$ \\
\hline
Uniform & $\frac{b-x}{b-a}$  & $\frac{(b-a)(as+a+b)}{(s+1)(s+2)}$\\

 Exponential & $\ e^{-\lambda x} $ & $\frac{1}{{s^2} {\lambda^2}}$  \\

 Pareto & $\frac{1}{x^{\lambda}}$ & $\frac{1}{{\lambda s}-2};{{\lambda s}>2}$  \\

 Lomax & $(\frac{m}{m+x})^n $&$\frac{m^2}{(2-ns)(1-ns)}$ \\

 Rayleigh & $\ e^{\frac{-x^2}{2{\sigma}^2}}$&$\frac{\sigma^2}{s}$ \\
\hline
\end{tabular}
\end{table}

Next, we study some properties of $C_s(W,F)$.   Suppose, for the random variable $X$ having distribution function $F$, let $Y=aX+b$ with $a >0$ and $b\geq 0$, then we can observe that $C_{s}(W,Y)$ is the sum of ${a^2}$ times WCREGF of $X$ and $ab$ times CREGF of $X$. That is
$$
C_{s}(W,Y)= {a^2}C_{s}(W,F)+ab C_{s}(X).
$$

The proportional hazards model  (Cox, 1972) estimates the effects of different covariates influencing the failure time of a system/subject. For more details on the theory and applications of the proportional hazards model, one may refer to Kalbfleisch and Prentice (2002). Let $X^*$ be a continuous non-negative random variable with survival function (proportional hazards model)
$$\bar{F}^*(t)= [\bar{F}(t)]^\theta,\,\theta>1,$$
where $\bar{F}^*$ is the survival function of $X^*$.
We can observe that the WCREGF of $X$ and $X^*$ are related by
$$C_s(W,F^*)=C_{s\theta}(W,F).$$

For a series system with $n$ independent and identical components each having lifetime $X_i,\,i=1,\ldots,n$, the lifetime of the system is given by $\min(X_1,\ldots,X_n)$. Denote $\bar F_s$ is the survival function of the random variable $\min(X_1,\ldots,X_n)$. Let $X_{(1)}$ be the first order statistic based on a random sample $X_1,...X_n$ from F. We have $\bar{F}_{s}(x)= [\bar{F}(x)]^n$, then $$C_s(W,F_s)=C_{
ns}(W,F).$$
In the following theorem, we obtain the relationship between $C_{s}(W,F)$ and Shannon's entropy measure.
\begin{theorem}
 Let $X$ be continuous non-negative random variable having Shannon's entropy $H(X)$ in (\ref{shannon}), then $$C_{s}(W,F)\geq \exp[H(X)+E(log X)-s].$$
\end{theorem}
\noindent{\bf Proof:}  By using log-sum inequality, we obtain
\begin{eqnarray*}
  \int_{0} ^ {\infty} f(x)\log\left(\frac{f(x)}{x (\bar{F}(x)  )^s}\right)dx &\geq & \int_{0}^{\infty}f(x) \log\left(\frac{\int_{0}^{\infty}f(x)dx} {\int_{0}^{\infty}x(\bar{F}(x))^sd x}\right) d x .\\
&\geq& -\log\left(\int_{0}^{\infty}x(\bar{F}(x))^sd x\right)
\\
&=&-\log (C_{s}(W,F)).
\end{eqnarray*}
From the above expression, we have the inequality stated in the theorem.
\section{Dynamic weighted cumulative residual entropy generating function}
In various fields such as reliability, survival analysis, economics, business, etc., the duration of a study period is a significant variable. In these instances, information-generating functions become dynamic, incorporating a time dependency. Inspired by the dynamic nature associated with time-dependent functions, this section introduces the dynamic version of WCREGF, which is known as the dynamic weighted cumulative residual entropy generating function (DWCREGF).  The DWCREGF is defined as
$$C_{s}(W,X; t)=\int_{t}^{\infty}x\left(\frac{\bar{F}(x)}{\bar{F}(t)}\right)^s dx,\,s>0.$$
\begin{remark}
 Clearly $ C_{s}(W,X; 0)= C_{s}(W,F)$, which is WCREGF.
\end{remark}

In Table 2, we give the expression for $C_{s}(W,F,t)$ for some well-known lifetime distributions.
\begin{table}[h]
  \centering
   \caption{Expression of $C_{s}(W,F,t)$ for some well-known distributions}
\begin{tabular}{cccccccc}
\hline\\
 Distribution & $\bar{F}(x)$ & $C_s(W,F,t)$ \\
\hline\\
Uniform & $\frac{b-x}{b-a}$ & $\frac{(b-t)[t(s+2)+(b-t)]}{(s+1)(s+2)}$\\

 Exponential & $\ e^{-\lambda x} $ & $\frac{(1+\lambda s t)}{( s \lambda)^2}$  \\

 Pareto & $\frac{1}{x^{\lambda}}$ & $\frac{t^2}{{\lambda s}-2};{{\lambda s}>2}$ \\

 Lomax & $(\frac{m}{m+x})^n $& $\frac{(m+t)}{(ns-2)}[t+\frac{(m+t)}{(ns-2)}]; ns-2>0$ \\

 Rayleigh & $\ e^{\frac{-x^2}{2{\sigma}^2}}$ & $\frac{\sigma^2}{s} $\\
\hline
\end{tabular}
\end{table}

\begin{definition}
Let $X$ be a non-negative random variable with survival function $\bar{F}(x)$, then the weighted mean residual lifetime (WMRL) is defined as
$$m_{F}(W;t)=\int_{t}^{\infty}x\frac{\bar{F}(x)}{\bar{F}(t)}dx.$$
\end{definition}
\noindent In particular
$m_{F}(W;0)=\int_{0}^{\infty}x\bar{F}(x)dx=\frac{1}{2}E(X^2).$
\begin{theorem}
Let $X$ be a non-negative random variable and let $C_{s}(W,F)<\infty,$ then
$$C_{s}(W,F)< m_{F}(W;0), s>0.$$\end{theorem}
\noindent{\bf Proof:}
The result follows from the definition of $C_{s}(W,F)$  and using the inequality $(\bar{F}(x))^s<\bar{F}(t)$,\,$t>0,  s>0$.
\begin{theorem}
For a non-negative random variable $X$ with survival function $\bar{F}(x)$, then
$$C_{s}(W,X;t)\leq m_{F}(W;t), s>0.$$\end{theorem}
 \noindent{\bf Proof:} For all $x\ge t$,  we have  $\bar{F}(x)\le \bar{F}(t)$. Hence using the definitions of $C_{s}(W,X;t)$ and $m_{F}(W;t)$, we obtain the desired result.

Next, we obtain the relationship between hazard rate and DWCREGF. By differentiating $C_{s}(W,X;t)$ with respect to $t$, we obtain the following result:
\begin{equation}\label{rent}
h(t)=\frac{t+C'_{s}(W,X;t)}{s.C_{s}(W,X;t)}.
\end{equation}
The following theorem shows that the dynamic weighted cumulative residual entropy generating function determines the distribution of $X$ uniquely.
\begin{theorem}
Let $X$ be non-negative random variable with density function $f(x)$, the survival function $\bar{F}(x)$ and the hazard rate $h(x)$. Then $C_{s}(W,X;t)$ uniquely determines the distribution of $X$.
\end{theorem}
\noindent{\bf Proof:}
Using the relationship between hazard rate and dynamic weighted cumulative residual entropy generating function in (\ref{rent}), we can write
\begin{equation*}
\frac{-d\log \bar{F}( t)}{dt}=\frac{t+C'_{s}(W,X;t)}{s.C_{s}(W,X;t)}.
\end{equation*}
Integrating for $t$  over the interval $(0,x)$, we have
$$\bar{F}( x)=\exp\left[{-\int_{0}^{x}\frac{t+C'_{s}(W,X;t)}{s.C_{s}(W,X;t)}}dt\right].$$
This shows that the knowledge of $ C_{s}(W,X;t) $ enables us to determine the distribution of $X$.
Now, suppose that $F(.)$ and $G(.)$ be two distribution functions such that
\begin{equation}\label{eq5}
C{s}(W,F;t)  = C_{s}(W,G;t).
\end{equation}
Differentiating both side of the above equation with respect to $t$, we obtain
$$-t+sh_{1}(t)C_{s}(W,F;t)= -t+sh_{2}(t)C_{s}(W,G;t),$$
where $h_{1}(t)$ and $h_{2}(t)$ are hazard rates corresponding to $F$ and $G$, respectively.\\
In view of equation (\ref{eq5}), we have $h_{1}(t)=h_{2}(t)$.
This implies that $C_{s}(W,F;t)$ determines the distribution of $X$ uniquely.

The next theorem shows that the dynamic weighted residual entropy generating function is independent of  $t$  if and only $X$ has Rayleigh distribution.
\begin{theorem}
  If $X$ is a non-negative continuous random variable having distribution function $F(x)$, then the dynamic weighted cumulative residual entropy generating function is independent of  $t$ if and only if $X$ has Rayleigh distribution.
\end{theorem}
\noindent{\bf Proof}. Let $C_{s}(W,F ; t)=k $, where $k$ is a positive constant. Then
$$
C_{s}^{\prime}(W,F ; t)=0 .
$$
From Equation (\ref{rent}), we obtain
$$
sk. h(t)=t.
$$
Or $$h(t)=\frac{t}{s k}=at,$$ where   $a$ is a constant.
And  $h(t)$ is the hazard rate corresponding to the Rayleigh distribution.

Conversely assume that, $X$ has  Rayleigh distribution  where the survival function is given by  $$\bar{F}(x)=\exp (\frac{-x^2}{2b^2}). $$ Therefore,
\begin{eqnarray*}
C_{s}(W,F ; t) &=&\frac{1}{e^\frac{-st^2}
{2b^2}} \int_{t}^{\infty} x{e^\frac{-sx^2}
{2b^2}} d x \\
&=&\frac{b^2}{s}=k,
\end{eqnarray*}
which is a constant. Accordingly, WDCREGF is independent of  $t$ if and only if $X$ has Rayleigh distribution.

\noindent We use this characterization property to develop goodness of fit test for Rayleigh distribution.
\begin{theorem}
If $m_{F}(W;t)$ is the weighted mean residual life of $X$, then the relation
\begin{equation}\label{rent5}
sC_{s}(W,X;t)=m_{F}(W;t)
\end{equation}
holds if and only if $X$ has the survival function $\bar{F}(t)=\exp (\frac{-t^2}{2b^2})$.
\end{theorem}
\noindent{Proof:}  Let  $X$ has Rayleigh distribution with survival function
 $$ \bar{F}(t)=\exp (\frac{-t^2}{2b^2}).$$
  Then $m_{F}(W;t)=b^2$ and $C_{s}(W,X;t)=\frac{b^2}{s}$ and the if part of the theorem holds.\\
  To prove the converse,  assume that equation (\ref{rent5}) holds.   Differentiating both sides of equation (\ref{rent5})  with respect to $t$, and using the definition of the hazard rate,  we obtain
  \begin{equation}\label{rent6}
  h(t).m_{F}(W;t)=t.
  \end{equation}
  But
  \begin{equation*}\label{rent2}
  \frac{d}{dt}m_{F}(W;t)=h(t).m_{F}(W;t)-t
  \end{equation*}
 Therefore, using the above two equations, we obtain
 $$\frac{d}{dt}m_{F}(W;t)=0 .$$
 Hence $m_{F}(W;t)$ is constant. Therefore equation (\ref{rent6}) becomes
 $$h(t)=\frac{t}{k}.$$
 This gives
 $$
 \bar{F}(x)=\exp(\frac{-t^2}{2k}),
 $$
  which is the survival function of Rayleigh distribution and hence the proof of the theorem.

\section{Test for Rayleigh distribution}\vspace{-0.1in}
In the previous section, we proved that the  WDCREGF characterizes the distribution of $X$.  Also, we proved that constant  WDCREGF is a characterization property of  Rayleigh distribution.   We develop a test for Rayleigh distribution against the decreasing DCREGF class.

Let $X_1,\ldots X_n$ be a random sample of size $n$ from $F$. We are interested in testing the null hypothesis\vspace{-0.1in}
$$H_0: X \text{ has Rayleigh distribution} $$
 against the alternative hypothesi\vspace{-0.1in}
 $$H_1: X \text{ has decreasing DCREGF and not Rayleigh.}$$
  For testing the above hypothesis first we define a departure measure that discriminates between null and alternative hypotheses.
  Note that $C_{s}(W, X ; t)$ is decreasing in $t$ if $C_{s}^{'}(W, X ; t)\le 0.$ That is,
  $$\frac{sf(t)}{\bar{F}^{s+1}(t)}\int_{t}^{\infty}x\bar{F}^s(x)dx-t\le 0$$
  or
  $${t\bar{F}^{s+1}(t)}-{sf(t)}\int_{t}^{\infty}x\bar{F}^s(x)dx\ge 0.$$
  Hence, we consider  a measure of departure $ \Delta(F)$ given by
\begin{eqnarray}\label{deltam}
\Delta(F)=\int_{0}^{\infty}\left({t\bar{F}^{s+1}(t)}-{sf(t)}\int_{t}^{\infty}x\bar{F}^s(x)dx\right)dt.
\end{eqnarray}Clearly, $\Delta(F)$ is zero under  $H_0$ and positive  under $H_1$. Accordingly,  $\Delta(F)$  can be considered as a measure of departure from $H_0$ towards  $H_1$. As the proposed test is based on U-statistics, first we express  $\Delta(F)$ in terms of expectation of the function of random variables. Observe that $\bar{F}^{n}(x)$ is the survival function of $\min(X_1,\ldots,X_n)$. For a non-negative random variable $E(\min(X_1,\ldots,X_{n}))=\int_{0}^{\infty}\bar{F}^n(x)dx$ and $E(\min(X_1,\ldots,X_{n})^2)=\int_{0}^{\infty}2x\bar{F}^n(x)dx.$
Consider
\begin{eqnarray}\label{delta}
\Delta(F)&=&\int_{0}^{\infty}\left({t\bar{F}^{s+1}(t)}-{sf(t)}\int_{t}^{\infty}x\bar{F}^s(x)dx\right)dt\nonumber\\
&=&\frac{1}{2}E(\min(X_1,\ldots,X_{s+1})^2)-\int_{0}^{\infty}{sf(t)}\int_{t}^{\infty}x\bar{F}^s(x)dxdt.
\end{eqnarray}
Changing the order of integration, from (\ref{delta}) we have
\begin{eqnarray}\label{delta1}
\Delta(F)&=&\frac{1}{2}E(\min (X_1,\ldots,X_{s+1})^2)-s\int_{0}^{\infty}x\bar{F}^s(x)\int_{0}^{x}f(t)dtdx\nonumber\\
&=&\frac{1}{2}E(\min(X_1,\ldots,X_{s+1})^2)-\frac{s}{2}\int_{0}^{\infty}x2\bar{F}^s(x)(1-\bar{F}(x))dx\nonumber\\\nonumber
&=&\frac{(s+1)}{2}E(\min(X_1,\ldots,X_{s+1}){^2})-\frac{s}{2}E(\min(X_1,\ldots,X_{s})^2).
\end{eqnarray}

We find the test statistic using the theory of U-statistics.  Consider a symmetric kernel
\begin{eqnarray*}
  h_1(X_1,\ldots,X_{s+1})= \frac{(s+1)}{2}\min(X_1,\ldots,X_{s+1})^2-\frac{1}{2(s+1)}\sum_{C_s}s\min(X_{i_{1}},\ldots,X_{i_{s}})^2,
\end{eqnarray*}where summation is over the set $C_s$ of all  combination of  $s$ integers $i_1<i_2<\ldots<i_s$ chosen from the set $(1,\ldots,s+1)$. Then $E(  h_1(X_1,\ldots,X_{s+1}))=\Delta(F)$.
Hence a U-statistic based  test statistic is given by
\begin{equation*}\label{test}
 \widehat{\Delta} =\frac{1}{C_{m,n}} \sum\limits_{C_{m,n}}h{(X_{i_1},X_{i_2},\cdots,X_{i_{s+1}})},
\end{equation*}where the summations is over the set  $C_{m,n}$ of all combinations of $(s+1)$ distinct elements $\lbrace i_1,i_2,\cdots,i_{s+1}\rbrace$ chosen from $\lbrace 1,2,\cdots,n\rbrace$.
We reject the null hypothesis $H_0$ against the alternative  $H_1$ for large value of $\widehat{\Delta}$. We obtain a critical region of the proposed test using the asymptotic distribution of $ \widehat{\Delta}$.

 In the next theorem, we state the asymptotic distribution of $\widehat{\Delta} $.
\begin{theorem}\vspace{-0.1in}
  As $n\rightarrow \infty$,  $\sqrt{n}(\widehat{\Delta}-\Delta(F))$ converges in distribution to normal random variable with mean zero and variance $(s+1)^2\sigma^2$, where $\sigma^2$ is given by
 \begin{eqnarray}\label{vart}
 \sigma^2&=&\frac{1}{4}Var\Big((s+1)X^2{{\bar F}^{s}}(X) + s(s+1)\int_0^X {{{y^2 \bar F}^{s -1}}(y)} d{F}(y)\nonumber\\
 &&\quad-\frac{s^2}{(s+1)}X^2{{\bar F}^{s-1}}(X) -\frac{(s-1)s^2}{(s+1)}\int_0^X {{{y^2 \bar F}^{s -2}}(y)} d{F}(y)\Big).
\end{eqnarray}
\end{theorem}
\noindent {\bf Proof:} By the central limit theorem for U-statistics,   asymptotic  distribution of $ \widehat{\Delta}$ is normal with mean $\Delta(F)$. The asymptotic variance is $(s+1)^2\sigma_1^2$ where $\sigma_1^2$ is given by (Lee, 2019)
\begin{equation}\label{var1}
\sigma_1^2= Var\left[E\left(h(X_{1},\ldots,X_{(s+1)})|X_{1}\right)\right].
\end{equation}Denote $Z=\min(X_2,X_3,...,X_s)$, then the distribution of $Z$ is given by $1-\bar{F}^{s-1}(x)$, where $\bar{F}(x)=1-F(x)$.  Consider
\begin{eqnarray*}\label{eq4}
E\left[ {\min \left( {x,{X_2},{X_3},...,{X_{s} }} \right)} \right]^2& =& E\left[ {x^2I(Z > x)} \right] + E\left[ {Z^2I(Z \le x)} \right]\nonumber\\&=&x^2{{\bar F}^{s- 1}}(x) + (s-1)\int_0^x {{{y^2 \bar F}^{s  - 2}}(y)} d{F}(y).
\end{eqnarray*}\vspace{-0.1in}Similarly, we have
\begin{eqnarray*}\label{eq4}
E\left[ {\min \left( {x,{X_2},{X_3},...,{X_{s+1} }} \right)} \right]^2&=&x^2{{\bar F}^{s}}(x) + s\int_0^x {{{y^2 \bar F}^{s  - 1}}(y)} d{F}(y).
\end{eqnarray*}\vspace{-0.1in}Hence \vspace{-0.1in}
\begin{eqnarray*}\label{eq4}
 &&\hskip-1.in E\left(h\big( {{X_1},{X_2},...,{X_s }}|X_1=x\big) \right)\\&=&(s+1)x^2{{\bar F}^{s}}(x) + s(s+1)\int_0^x {{{y^2 \bar F}^{s -1}}(y)} d{F}(y)\\
 &&-\frac{s^2}{(s+1)}x^2{{\bar F}^{s-1}}(x) -\frac{(s-1)s^2}{(s+1)}\int_0^x {{{y^2 \bar F}^{s -2}}(y)} d{F}(y)+\frac{sk}{2(s+1)},
\end{eqnarray*}where $k=E[\min(X_2,\ldots, X_{s+1})]^2$, a constant.
\noindent Therefore, from (\ref{var1}) we obtain the variance expression specified in the statement of the theorem.

Under the null hypothesis $H_0$, $\Delta{(F)}=0$. Hence we have the following corollary.
\begin{corollary}
 Under $H_0$, as $n\rightarrow \infty$,  $\sqrt{n}\widehat{\Delta}$ converges in distribution to a Gaussian random variable with mean zero and variance $\sigma_0^2$, where $\sigma_0^2$ is the value of the $\sigma^2$ evaluated under $H_0$.
\end{corollary}
An asymptotic critical region of the test can be obtained using Corollary 1. Let  $\widehat\sigma_0^2$ be a consistent estimator of the asymptotic null variance $\sigma_0^2$. We reject the null hypothesis $H_{0}$ against the alternative hypothesis $H_{1}$ at an approximate significance level $\alpha$, if
\begin{equation*}
 \frac{ \sqrt{n}\widehat{\Delta}}{{\widehat\sigma_0}}>Z_{\alpha},
  \end{equation*}
where $Z_{\alpha}$ is the upper $\alpha$-percentile point of the standard normal distribution.

\section{Simulation Studies}
We conduct an extensive Monte Carlo simulation study to assess the finite sample performance of the proposed test using the empirical type I error and power of the test. The simulation is done using the R software and repeated ten thousand times using different sample sizes.

The empirical type I error and power of test are compared with the classical goodness of fit tests including the Anderson Darling (AD) test, Cramer von Mises (CvM) test and  Kolmogrv-Sminrnov (KS) test and other goodness of fit tests, specifically developed for Rayleigh distribution.
The proposed test can be performed with different choices of the parameter $s.$ In simulation studies, we observe that the performance of the proposed test is optimum for the choice $s=1.$ Hence for the comparison purpose with other competent tests, we choose $s=1.$
We generate lifetimes of different sample sizes, ($n=10,20,30,40,50$) to calculate the empirical type I error and power of the proposed test and the other competent tests.

Finding a  consistent estimator of the null variance $\sigma_{0}^2$ is difficult. Hence we find the critical region of the proposed test using Monte Carlo simulation. We determine $\alpha$-th  quantiles $c_1$ in such a way that $P(\widehat\Delta>c_1)=\alpha$.
The algorithm used to find the empirical power of the proposed test is summarised as follows.
\begin{enumerate}
\item Generate lifetime data from the desired alternative and calculate the value of $\widehat\Delta$.
\item Generate lifetime data from the Rayleigh distribution and calculate the value of $\widehat\Delta$.
\item Repeat Step 2 10000 times and determine the critical point.
\item Repeat Steps 1-3 10000 times and calculate empirical power as the proportion of significant test statistics.

\end{enumerate}

\begin{table}[h]
\caption{ Comparison of empirical power for different alternatives ($\alpha=0.01$) }
\label{Tab:2}
\begin{tabular}{lrrrrrrrrr} \hline \noalign{\smallskip}
& $n$ & SSS & KS & CvM & AD & AH & MI & JH & VTS\\
\noalign{\smallskip} \hline \noalign{\smallskip}

Rayleigh(1) & 10 & 0.0121 & 0.0175 & 0.0193 & 0.0119 & 0.0114 & 0.0118 & 0.0086 & 0.0115 \\
& 20 &   0.0115 & 0.0142 & 0.0161 & 0.0110 & 0.0109 & 0.0112 & 0.0090 & 0.0111 \\
& 30 &  0.0110 & 0.0128& 0.0127 & 0.0105 & 0.0092 & 0.0107 & 0.0108 & 0.0106 \\
& 40 &   0.0094 & 0.0116 & 0.0119 & 0.0095 & 0.0095 & 0.0094 & 0.0105 & 0.0095 \\
& 50 &   0.0096 & 0.0111 & 0.0114 & 0.0102 & 0.0104 & 0.0096 & 0.0105 & 0.0098 \\

\noalign{\smallskip} \hline \noalign{\smallskip}
Weibull (2, 2) & 10 & 0.3589 & 0.2456 &0.1892 &0.5673 &0.1298& 0.3974& 0.1988& 0.2387 \\
& 20 &0.6586 & 0.5784& 0.5181& 0.7172& 0.3912& 0.6313 &0.2987& 0.5218 \\
& 30 & 0.8287 & 0.8241 &0.8293& 0.8518& 0.4394& 0.7052 &0.4688& 0.6959 \\
& 40 &0.9461 & 0.9396& 0.9422& 0.9219& 0.5865 &0.7366 &0.6551 &0.8277 \\
& 50 & 0.9869& 0.9797& 0.9835 &0.9723& 0.7893 &0.8248 &0.8102& 0.8978  \\

\noalign{\smallskip} \hline \noalign{\smallskip}
Pareto (2,2) & 10 & 0.9237 &  0.9704 & 1.0000 &0.9968& 0.2714& 0.5288& 0.2981 & 0.5593 \\
& 20 & 0.9692 & 1.0000 & 1.0000 & 1.0000 & 0.3892 &0.5620& 0.6118& 0.7785    \\
& 30 & 0.9814 & 1.0000 & 1.0000& 1.0000 &0.5998 &0.8229& 0.8367 &0.8387 \\
& 40 & 0.9843& 1.0000 & 1.0000& 1.0000 &0.7461 &0.8736& 0.8848& 0.9166 \\
& 50 & 0.9885& 1.0000 & 1.0000 & 1.0000 &0.8823 &0.9159& 0.9285& 0.9459  \\

\noalign{\smallskip} \hline \noalign{\smallskip}
Log-normal (1,1) & 10 &0.6927 &0.9755 &0.9812& 0.9854 &0.3414& 0.3563 &0.2856 &0.5470  \\
& 20  & 0.9418 & 1.0000& 1.0000 &1.0000 &0.6534 &0.6386& 0.4196& 0.8733  \\
& 30 & 0.9816 & 1.0000 &1.0000& 1.0000 &0.8275& 0.8042 &0.6954 & 0.9671  \\
& 40 & 0.9924 & 1.0000 &1.0000 &1.0000 &0.9266& 0.9151& 0.8566 &0.9935 \\
& 50 & 0.9962 & 1.0000& 1.0000 &1.0000 &0.9584& 0.9554 &0.9532 &0.9985   \\

\noalign{\smallskip} \hline \noalign{\smallskip}
Half-normal (0.5) & 10 & 0.3727 &  0.5036 & 0.6457 & 0.5854 & 0.1987 & 0.2795 & 0.1643 & 0.3654 \\
& 20 &  0.6463 & 0.8925& 0.9382 &0.8361& 0.3936& 0.3317 &0.2666& 0.5646 \\
& 30 & 0.8241& 0.9831 &0.9934& 0.9453 &0.5404 &0.4218& 0.4097 &0.7438  \\
& 40 & 0.9102&  0.9980& 0.9994 & 0.9842& 0.7415 & 0.5762 & 0.5924& 0.8589  \\
& 50 & 0.9598 &0.9995 &1.0000 &0.9928& 0.7845& 0.6517 &0.7215 & 0.9283 \\

\noalign{\smallskip} \hline \noalign{\smallskip}
LFR (0.5) & 10 & 0.8779 & 0.7316& 0.8175& 0.7597 &0.1911& 0.1939 &0.2841& 0.7658   \\
& 20 & 0.9956 &  0.9874& 0.9952& 0.9635 & 0.2625 & 0.3543 & 0.3184& 0.8279\\
& 30 & 1.0000 & 1.0000& 1.0000 &0.9898 &0.5421& 0.4145 &0.6478& 0.8811 \\
& 40 & 1.0000 &1.0000& 1.0000& 1.0000 &0.7132 &0.6545& 0.7933 &0.9221
\\
& 50 & 1.0000 &1.0000& 1.0000& 1.0000& 0.7818& 0.7561& 0.8882& 0.9512
 \\
\hline
\end{tabular}
\end{table}

\begin{table}[h]
\caption{Comparison of empirical power for different alternatives ($\alpha=0.05$)}
\label{Tab:3}
\begin{tabular}{lrrrrrrrrr} \hline \noalign{\smallskip}
& $n$ & SSS & KS & CvM & AD & AH & MI & JH & VTS\\
\noalign{\smallskip} \hline \noalign{\smallskip}
Rayleigh (1) & 10 & 0.0479 & 0.0563 & 0.0575 & 0.0529 & 0.0442 & 0.0525 & 0.0521 & 0.0532    \\
& 20 & 0.0488 &0.0551 & 0.0561 & 0.0518 & 0.0471 & 0.0518 & 0.0516 & 0.0511  \\
& 30 & 0.0512 &0.0526 & 0.0531 & 0.0478 & 0.0482 & 0.0487 & 0.0515 & 0.0509   \\
& 40 & 0.0509& 0.0512 & 0.0519 & 0.0489 & 0.0484 & 0.0492 & 0.0478 & 0.0489   \\
& 50 & 0.0503 &0.0491 & 0.0489 & 0.0490 & 0.0489 & 0.0505 & 0.0494 & 0.0495 \\
\noalign{\smallskip} \hline \noalign{\smallskip}
Weibull (2, 2) & 10 & 0.5758 &0.3204 &0.2616 &0.6567 &0.2041& 0.4197& 0.2292 &0.3528  \\
& 20 &0.8399& 0.6178& 0.5718& 0.7728& 0.4109& 0.6731 &0.3206& 0.6521  \\
& 30 & 0.9393& 0.8424 &0.8629& 0.8846& 0.5199& 0.7405& 0.5368 &0.7695 \\
& 40 & 0.9771 &0.9539& 0.9642& 0.9423& 0.286& 0.7836 &0.6904& 0.8921  \\
& 50 & 0.9935 &0.9979 &0.9983& 0.9972& 0.8008& 0.8824& 0.8567& 0.9487   \\

\noalign{\smallskip} \hline \noalign{\smallskip}
Pareto (2, 2) & 10 & 0.9395& 1.0000 &1.0000 &1.0000& 0.3109& 0.6351& 0.3326 &0.6571  \\
& 20 & 0.9704 & 1.0000 &1.0000 &1.0000 &0.4257& 0.7512 &0.6596 &0.8043  \\
& 30 & 0.9904 &1.0000& 1.0000& 1.0000& 0.6323& 0.7599& 0.8672 &0.8976 \\
& 40 & 0.9984 & 1.0000& 1.0000& 1.0000 &0.8256& 0.8226& 0.9098& 0.9419  \\
& 50 &1.0000 &1.0000 &1.0000 &1.0000 &0.9382& 0.8641& 0.9695& 0.9763  \\

\noalign{\smallskip} \hline \noalign{\smallskip}
Log-normal (1,1) & 10 &  0.7438& 0.9934& 0.9947 &0.9951& 0.5675 &0.4399 &0.3934 &0.7367  \\
& 20 &  0.9511 & 1.0000 &1.0000& 1.0000 &0.7763& 0.6845 &0.6445& 0.9261 \\
& 30 & 0.9891& 1.0000 &1.0000 &1.0000& 0.8490& 0.8516& 0.8218 &0.9853  \\
& 40 & 0.9972 & 1.0000& 1.0000 &1.0000& 0.9365& 0.9494& 0.9051 &0.9976
 \\
& 50 & 1.0000 &1.0000& 1.0000&1.0000 &0.9876 & 0.9794 &0.9739 &1.0000
 \\

\noalign{\smallskip} \hline \noalign{\smallskip}
Half-normal (0.5) & 10 &  0.5075& 0.7254 &0.7863& 0.7867 &0.3291 &0.3214& 0.2804 &0.4712  \\
& 20 & 0.7888 & 0.9569& 0.9684& 0.9608& 0.4505& 0.4379 &0.3898& 0.6529  \\
& 30 &0.9073 & 0.9971& 0.9988& 0.9895& 0.6446 &0.5355 &0.5378& 0.7986  \\
& 40 & 0.9545 & 0.9998 &1.0000 &0.9998& 0.7935& 0.6472 &0.6489& 0.8921    \\
& 50 & 0.9855 &  1.0000& 1.0000 &1.0000& 0.8528 & 0.7828 &0.7942 &0.9883
\\

\noalign{\smallskip} \hline \noalign{\smallskip}
LFR (0.5) & 10 & 0.9521& 0.9097& 0.9373& 0.9584 & 0.2871 &0.2448 & 0.3109 &0.8201

\\
& 20 & 0.9998 & 0.9958 &0.9982& 0.9972 & 0.4291 &0.5348 &0.5621 &0.8827

 \\
& 30 & 1.0000& 1.0000& 1.0000& 1.0000 & 0.6532& 0.6689&0.7818&0.9237

 \\
& 40 &  1.0000& 1.0000& 1.0000& 1.0000&0.7881 &0.7387&0.8356&0.9459
\\
& 50 & 1.0000& 1.0000 & 1.0000 & 1.0000&0.8289 &0.8636&0.8902&0.9752\\
\hline
\end{tabular}
\end{table}
To find the empirical power, lifetime random variables are generated from different choices of alternatives including Weibull, log-normal, Pareto, half-normal and linear failure rate (LFR) distributions, where the distribution functions are:
    \begin{itemize}
        \item Weibull distribution: $F_{1}(x)=1-e^{(-x/\lambda)^{k}}$, $x>0$, $k,\,\lambda>0$.
          \item  Pareto distributions: $F_{2}(x)=(\lambda/x)^{\alpha}$ $x>0$, $\alpha,\, \lambda>0$.
    \item  Lognormal distribution: $F_{3}(x)=\Phi(\frac{\ln x-\mu}{\sigma})$, $x>0$, $-\infty<\mu<\infty,\, \sigma^2>0$
    where $\Phi(x)$ is the cumulative distribution function of the standard normal random variable.
    \item Half-normal distribution: $F_{4}(x)= erf(\frac{x}{\sigma\sqrt{2}})$, $\sigma >0$, where \textit{erf} is the error function.
    \item LFR distribution: $F_{5}(x)= 1- \exp(-\frac{-\lambda x^2}{2})$~~~ $x \ge 0$,~~~~$\lambda>0$.
    \end{itemize}
We compare the performance of our test (denoted by SSS)  with goodness of fit test for Rayleigh distribution proposed by Ahrari et al. (2022) (AH), Meintanis and Illipoulos (2003) (MI), Jahanshahi et al. (2016) (JH) and Vaishakh et al. (2022) (VTS) as well as with the classical tests  Kolmogorov Smirnov (KS), Cramer von Mises (CvM) and Anderson-Darling (AD) test.

We report the results of the empirical power comparison in Tables 1 and 2. Table 1 gives the power of $\widehat \Delta$ for different alternatives at significance level $\alpha=0.01$ and Table 2 presents the same when $\alpha=0.05$.  It is evident from Tables 1 and 2 that, the empirical type I error of the test approaches the chosen significance level as sample size increases.  Also, the proposed test exhibits good power against all choices of alternatives which increases with sample size. The results from the simulation study ensure that the proposed test is superior to other competent tests.

\vspace{-0.3in}
\section{Data Analysis}
We illustrate the proposed testing procedure with two real data sets. To find the p-values, we use the following parametric bootstrap procedure.
\begin{enumerate}
\item Obtain the moment estimate of the parameter of Rayleigh distribution from the observed data.
\item  Generate 10000 bootstrap samples from Rayleigh distribution having parameter estimated in Step 1.
  \item Calculate the test static value for each sample generated in the previous step.
            \item Calculate the bootstrap p-value as the proportion of the test static value calculated in Step 3 that is greater than the test static estimated for the original data.
        \end{enumerate}

{\bf Data 1 }: We analyze the data, given in Caroni (2002), which represents the failure times of 23 ball bearings, to illustrate the proposed test. These failure times are: {\textit{17.88, 28.92, 33.00, 41.52, 42.12, 45.60, 48.48, 51.84, 51.96, 54.12, 55.56, 67.80, 67.80, 67.80, 68.88, 84.12, 93.12, 98.64, 105.12, 105.84, 127.92, 128.04, 173.40.}}  This data set originally used by Lieblein and Zelen (1956), is also discussed in Caroni (2002). The observations are the number of million revolutions before failure for each 23 ball bearings ordered according to life endurance. Later, the data is analysed by many authors including Kim and Han (2009), Dey and Dey (2014) and Vaisakh et al. (2023) in studying Rayleigh distribution. The bootstrap p-value, using the above-described computational algorithm is obtained as 0.6425. Accordingly, we can conclude that,  the hypothesis that the data follows Rayleigh distribution can be accepted for this data.  We can note that this result agrees with the results from previous studies on Rayleigh distribution.

{\bf Data 2 }: To illustrate the performance of the proposed test, we consider the data on survival times in weeks for 20 male rats that were exposed to a high level of radiation. The data are due to Furth et al. (1959)
and have been studied by Lawless (2003).  The complete data is given in Example 4.2.1 (Page 168) of Lawless (2003). The bootstrap p-value is obtained as 0.0152. Hence we reject the null hypothesis that the data follows Rayleigh distribution at a 5\% level of significance. Vaisakh et al. (2023) showed that survival time data follow a gamma distribution.
\section{Conclusion}
The generating function approach to entropy become popular in recent times. In this paper, we defined the weighted cumulative residual entropy generating function and studied its properties. The dynamic version of the proposed generating function is also given. We proved that the dynamic weighted cumulative residual entropy generating function (DWCREGF) determines the distribution uniquely. We proved that the constant DWCREGF characterizes the Rayleigh distribution. Using this characterization result, we developed a goodness of fit test for Rayleigh distribution and studied its properties. A Monte Carlo simulation study showed that the proposed test has a well-controlled error rate and has good power for various alternatives.  The proposed test is illustrated using two real data sets.


\begin{thebibliography}{999}
 \bibitem{} Ahrari, V., Baratpour, S., Habibirad, A. and Fakoor, V. (2022). Goodness of fit tests for Rayleigh distribution based on quantiles.{\em Communications in Statistics-Simulation and Computation}, 51,
341--357.
 \bibitem{} Balakrishnan, N., Buono, F. and Longobardi, M. (2022). On weighted extropies. {\em Communications in Statistics-Theory and Methods}, 51,6250--6267.
 \bibitem{}
Caroni, C. (2002). The correct ball bearings data. {\em Lifetime Data Analysis},  8, 395--399.
 \bibitem{} Capaldo, M., Di Crescenzo, A., and  Meoli, A. (2023). Cumulative information generating function and generalized Gini functions. {\em Metrika,} 1--29.
  \bibitem{} Chakraborty, S. and Pradhan, B. (2023). On weighted cumulative Tsallis residual and past entropy measures. {\em Communications in Statistics-Simulation and Computation}, 52, 2058--2072.
 \bibitem{} Clark, D. E. (2020). Local entropy statistics for point processes. {\em IEEE Transactions on Information Theory}, 66, 1155--1163.
  \bibitem{} Cox, D. R. (1972). Regression models and life‐tables. {\em Journal of the Royal Statistical Society: Series B}, 34, 187--202.

  \bibitem{}
        Dey, S. and Dey, T. (2014). Statistical inference for the Rayleigh distribution under progressively Type-II censoring with binomial removal. {\em Applied Mathematical Modelling},  38, 974--982.

\bibitem{}Di Crescenzo, A. and Longobardi, M. (2009). On cumulative entropies. {\em Journal of Statistical Planning and Inference}, 139, 4072--4087.
\bibitem{}
    Furth, J., Upton, A. C.,  and  Kimball, A. W. (1959). Late pathologic effects of atomic detonation and their pathogenesis. {\em Radiation Research Supplement} , 1, 243--264

\bibitem{} 	Golomb, S. (1966). The information generating function of probability distribution. {\em IEEE Transactions on Information Theory}, 12, 75--79.
 \bibitem{} Jahanshahi, S. M. A., Rad, A. H. and Fakoor, V. (2016). A goodness-of-fit test for Rayleigh distribution based on Hellinger distance. {\em Annals of Data Science}, 3, 401--411.
 \bibitem{} Kalbfleisch, J. D. and Prentice, R. L. (2011). {\em The Statistical Analysis of Failure Time Data}. John Wiley \& Sons, New Jersey.
\bibitem{}Kharazmi, O.  and Balakrishnan, N. (2021). Cumulative residual and relative cumulative residual Fisher information and their properties. {\em IEEE Transactions on Information Theory}, 67, 6306--6312.
\bibitem{} Kharazmi, O. and Balakrishnan, N. (2021). Jensen-information generating function and its connections to some well-known information measures. {\em Statistics \& Probability Letters}, 170, 108995.
\bibitem{} Kharazmi, O. and Balakrishnan, N. (2021). Cumulative and relative cumulative residual information generating measures and associated properties. {\em Communication in Statistics-Theory and Methods}, 5260--5273.
\bibitem{}
    Kim, C. and Han, K. (2009). Estimation of the scale parameter of the Rayleigh distribution with multiply type–II censored sample. {\em Journal of Statistical Computation and Simulation},  79, 965--976.
\bibitem{lee2019u} Lee, A. J. (2019). {\em U-statistics: Theory and Practice},  Routledge, New York.
 \bibitem{}
       Lieblein, J. and  Zelen, M. (1956). Statistical investigation of the fatigue life of deep-groove ball bearings. { \em Journal of Research of the National Bureau of Standards},  57, 273--316.

 \bibitem{}  Meintanis, S. and Iliopoulos, G. (2003). Tests of fit for the Rayleigh distribution based on the
empirical Laplace transform. {\em Annals of the Institute of Statistical Mathematics}, 55, 137--151.

\bibitem{} Mirali, M., Baratpour, S. and Fakoor, V. (2016). On weighted cumulative residual entropy. {\em Communications in Statistics-Theory and Methods}, 46, 2857--2869.
\bibitem{}Mirali, M. and Baratpour, S (2017). Some results on weighted cumulative entropy. {\em Journal of the Iranian Statistical Society}, 16, 21--32.
\bibitem{}	Rao, M.,  Chen, Y.,  Vemuri, B.,  Wang, F. (2004). Cumulative residual entropy: A new measure of information. {\em IEEE Transactions on Information  Theory}, 50, 1220--1228.
 \bibitem{} Saha, S. and Kayal, S. (2023). General weighted information and relative information generating functions with properties. arXiv preprint arXiv:2305.18746.
\bibitem{} Shannon, C. E. (1948). A mathematical theory of communication. {\em The Bell System Technical Journal}, 27, 379--423.
    \bibitem{} Smitha S., Sudheesh, K. K.,
 and Sreedevi, E. P. (2023). Dynamic cumulative residual entropy generating function and its properties. {\em Communications in Statistics-Theory and Methods}, 1--26.
\bibitem{}
Sudheesh, K. K.,  Sreedevi, E. P. and  Balakrishnan, N. (2022). A generalized measure of cumulative residual entropy. {\em Entropy}, 24, 444.
\bibitem{} Vaisakh, K. M., Xavier, T. and Sreedevi, E. P. (2023). Goodness of fit test for Rayleigh distribution with censored observations. {\em Journal of the Korean Statistical Society}, 52, 794--815.
\bibitem{} Zamani, Z., Kharazmi, O. and Balakrishnan, N. (2022). Information generating function of record values. {\em
Mathematical Methods of Statistics}, 31, 120--133.





\end{thebibliography}
\end{document}